\documentclass[a4paper,12pt]{article}
\usepackage[T1]{fontenc}
\usepackage{geometry}
\usepackage{textcomp}
\usepackage{amsmath}
\usepackage{graphicx,epsfig,color}
\usepackage{amsmath,amssymb,amsthm}
\usepackage{fancyhdr}
\DeclareMathAlphabet{\mathcalligra}{T1}{calligra}{m}{n}

\makeatletter

\newtheorem{thm}{Theorem}
\newtheorem{lem}[thm]{Lemma}
\newtheorem{cor}[thm]{Corollary}
\newtheorem{prop}[thm]{Proposition}
\theoremstyle{remark}
\newtheorem{remark}[thm]{Remark}
\DeclareMathOperator{\cg}{\textbf{[}}
\DeclareMathOperator{\cd}{\textbf{]}}

\newcommand{\cro}[1]{\cg {#1} \cd}

\def\N{\ensuremath{\mathbb{N}}}
\def\P{\ensuremath{\mathbb{P}}}

\newcommand{\reff}[1]{(\ref{#1})}

\def\N{\ensuremath{\mathbb{N}}}
\def\P{\ensuremath{\mathbb{P}}}

\date{\today}
\title{ \bf Stable continuous-state branching processes with immigration and Beta-Fleming-Viot processes with immigration.}
\author{ \bf Cl\'ement Foucart \\
\emph {Laboratoire de Probabilit\'es et Mod\`eles Al\'eatoires} \\
\emph {Universit\'e Pierre et Marie Curie}\\ 
\emph{4 Place Jussieu- 75252 Paris Cedex 05- France}
\\
\and 
\bf Olivier H\'enard\\
\emph {Universit\'e Paris-Est, CERMICS}\\ 
\emph{6 et 8 avenue Blaise Pascal}\\
\emph{77455 Marne-la-Vall\'ee - France}
}
\begin{document}
\maketitle{}
\begin{center}
DRAFT VERSION
\end{center}
\begin{abstract}
Branching processes and Fleming-Viot processes are two main models in stochastic population theory. Incorporating an immigration in both models, we generalize the results of Shiga (1990) and Birkner \textit{et al.} (2005) which respectively connect the Feller diffusion with the classical Fleming-Viot process and the $\alpha$-stable continuous state branching process with the $Beta(2-\alpha,\alpha)$-generalized Fleming-Viot process. In a recent work, a new class of probability-measure valued processes, called $M$-generalized Fleming-Viot processes with immigration, has been set up in duality with the so-called $M$-coalescents. The purpose of this article is to investigate the links between this new class of processes and the continuous-state branching processes with immigration. In the specific case of the $\alpha$-stable branching process conditioned to be never extinct, we get that its genealogy is given, up to a random time change, by a $Beta(2-\alpha, \alpha-1)$-coalescent. 
\end{abstract}
 \vspace{9pt} \noindent {\bf Key words.}
{Measure-valued processes}, {Continuous-state branching processes}, {Fleming-Viot processes}, {Immigration}, {Beta-Coalescent}, {Generators}, {Random time change}.
\par \vspace{9pt}
\par \vspace{9pt} \noindent {\bf e-mails.} {clement.foucart@etu.upmc.fr, henardo@cermics.enpc.fr}
\newpage
\section{Introduction} 
The connections between the Fleming-Viot processes and the continuous-state branching processes have been intensively studied. Shiga established in 1990 that a Fleming-Viot process may be recovered from the \textit{ratio process} associated with a Feller diffusion up to a random time change, see \cite{Shiga}. This result has been generalized in 2005 by Birkner \textit{et al} in \cite{Birk} in the setting of $\Lambda$-generalized Fleming-Viot processes and continuous-state branching processes (CBs for short). In that paper they proved that the ratio process associated with an $\alpha$-stable branching process is a time-changed $Beta(2-\alpha,\alpha)$-Fleming-Viot process for $\alpha \in (0,2)$. The main goal of this article is to study such connections when immigration is incorporated in the underlying population. The continuous-state branching processes with immigration (CBIs for short) are a class of time-homogeneous Markov processes with values in $\mathbb{R}_{+}$. They have been introduced by Kawazu and Watanabe in 1971, see \cite{KAW}, as limits of rescaled Galton-Watson processes with immigration. These processes are characterized by two functions $\Phi$ and $\Psi$ respectively called the immigration mechanism and the branching mechanism. A new class of measure-valued processes with immigration has been recently set up in 
Foucart \cite{coaldist}. These processes, called $M$-generalized Fleming-Viot processes with immigration ($M$-GFVIs for short) are valued in the space of probability measures on $[0,1]$. The notation $M$ stands for a couple of finite measures $(\Lambda_{0}, \Lambda_{1})$ encoding respectively the rates of immigration and of reproduction. The genealogies of the $M$-GFVIs are given by the so-called $M$-coalescents. These processes are valued in the space of the partitions of $\mathbb{Z}_+$, denoted by $\mathcal{P}^{0}_{\infty}$.
\\

In the same manner as Birkner \textit{et al.} in \cite{Birk}, Perkins in \cite{Perkins} and Shiga in \cite{Shiga}, we shall establish some relations between continuous-state branching processes with immigration and $M$-GFVIs. A notion of continuous population with immigration may be defined using a flow of CBIs in the same spirit as Bertoin and Le Gall in \cite{LGB0}. This allows us to compare the two notions of continuous populations provided respectively by the CBIs and by the $M$-GFVIs. Using calculations of generators, we show 
in Theorem \ref{thm1} that the following self-similar CBIs admit time-changed $M$-GFVIs for ratio processes: 
\begin{itemize} 
\item the Feller branching diffusion with branching rate $\sigma^{2}$ and immigration rate $\beta$ (namely the CBI with $\Phi(q)=\beta q$ and $\Psi(q)=\frac{1}{2}\sigma^{2}q^{2}$) which has for ratio process a time-changed $M$-Fleming-Viot process with immigration 
where $M=(\beta\delta_{0}, \sigma^{2}\delta_{0})$, 
\item the CBI process with $\Phi(q)=d'\alpha q^{\alpha-1}$ and $\Psi(q)=dq^{\alpha}$ for some $d,d'\geq 0$, $\alpha \in (1,2)$ which has for ratio process a time-changed $M$-generalized Fleming-Viot process with immigration 
where $M=\left(c'Beta(2-\alpha,\alpha-1), c Beta(2-\alpha,\alpha)\right)$, $c'=\frac{\alpha(\alpha-1)}{\Gamma(2-\alpha)}d'$ and $c=\frac{\alpha(\alpha-1)}{\Gamma(2-\alpha)}d$.
\end{itemize}

We stress that the CBIs may reach $0$, 
see Proposition \ref{extCBI}, in which case the $M$-GFVIs involved describe the ratio process up to this hitting time only.
When $d=d'$ or $\beta=\sigma^2$, the corresponding CBIs are respectively the $\alpha$-stable branching process and the Feller branching diffusion \textit{conditioned to be never extinct}. In that case, the $M$-coalescents are genuine $\Lambda$-coalescent viewed on $\mathcal{P}^{0}_{\infty}$. We get respectively a $Beta(2-\alpha, \alpha-1)$-coalescent when $\alpha \in (1,2)$ and a Kingman's coalescent for $\alpha=2$, see 
Theorem \ref{coalescent}. This differs from the $\alpha$-stable branching process \textit{without immigration} (already studied in \cite{Birk}) for which the coalescent involved is a $Beta(2-\alpha, \alpha)$-coalescent.
\\
Last, ideas provided to establish our main theorem have been used by Handa \cite{Handa} to study stationary distributions for another class of generalized Fleming-Viot processes. 
\\
\\
\textbf{Outline.} The paper is organized as follows. In Section 2, we recall the definition of a continuous-state branching process with immigration and of an $M$-generalized Fleming-Viot process with immigration. We describe briefly how to define from a flow of CBIs a continuous population represented by a measure-valued process. We state in Section 3 the connections between the CBIs and $M$-GFVIs, mentioned in the Introduction, and study the random time change. Recalling the definition of an $M$-coalescent, we focus in Section 4 on the genealogy of the $M$-GFVIs involved. We establish that, when the CBIs correspond with CB-processes conditioned to be never extinct, the $M$-coalescents involved are actually classical $\Lambda$-coalescents. We identify them and, as mentioned, the $Beta(2-\alpha,\alpha-1)$-coalescent arises. In Section 5, we compare the generators of the $M$-GFVI and CBI processes and 
prove the main result.
\section{A continuous population embedded in a flow of CBIs and the $M$-generalized Fleming-Viot with immigration}
\subsection{Background on continuous state branching processes with immigration} \label{CBI}
We will focus on critical continuous-state branching processes with immigration characterized by two functions of the variable $q\geq 0$:
\begin{align*}
\Psi(q)&=\frac{1}{2}\sigma^{2}q^{2}+\int_{0}^{\infty}(e^{-qu}-1+qu)\hat{\nu_{1}}(du)\\
\Phi(q)&=\beta q+\int_{0}^{\infty}(1-e^{-qu})\hat{\nu_{0}}(du)
\end{align*}
where $\sigma^{2}, \beta \geq 0$ and $\hat{\nu_{0}}$, $\hat{\nu_{1}}$ are two L\'evy measures such that $\int_{0}^{\infty} (1\wedge u) \hat{\nu_{0}}(du)<\infty$ and $\int_{0}^{\infty} (u\wedge u^{2}) \hat{\nu_{1}}(du)<\infty$. The measure $\hat{\nu_{1}}$ is the L\'evy measure of a spectrally positive L\'evy process which characterizes the reproduction. The measure $\hat{\nu_{0}}$ characterizes the jumps of the subordinator that describes the arrival of immigrants in the population. The non-negative constants $\sigma^2$ and $\beta$ correspond respectively to the continuous reproduction and the continuous immigration. Let $\mathbb{P}_{x}$ be the law of a CBI $(Y_{t}, t\geq 0)$ started at $x$, and denote by $\mathbb{E}_{x}$ the associated expectation. The law of 
the Markov process $(Y_{t},t\geq 0)$ can then be characterized by the Laplace transform of its marginal as follows: for every $q>0$ and $x\in \mathbb{R}_{+}$, $$\mathbb{E}_{x}[e^{-qY_{t}}]=\exp \left(-xv_{t}(q)-\int_{0}^{t}\Phi(v_{s}(q))ds \right)$$ 
where $v$ is the unique non-negative solution of $\frac{\partial}{\partial t}v_{t}(q)=-\Psi(v_{t}(q))$, $v_{0}(q)=q$.\\
\\
The pair $(\Psi, \Phi)$ is known as the branching-immigration mechanism. A CBI process $(Y_{t}, t\geq 0)$ is said to be conservative if for every $t>0$ and $x\in[0,\infty[, \mathbb{P}_{x}[Y_{t}<\infty]=1$. A result of Kawazu and Watanabe \cite{KAW} states that $(Y_{t}, t\geq 0)$ is conservative if and only if for every $\epsilon>0$ $$\int_{0}^{\epsilon}\frac{1}{|\Psi(q)|}dq=\infty.$$ 
Moreover, we shall say that the CBI process is \textit{critical} when $\Psi'(0)=0$: 
in that case, the CBI process is necessarily conservative.
We follow the seminal idea of Bertoin and Le Gall in \cite{LGB0} to define a genuine continuous population model with immigration on $[0,1]$ associated with a CBI. Emphasizing the r\^ole of the initial value, we denote by $(Y_{t}(x), t\geq 0)$ a CBI started at $x\in \mathbb{R}_{+}$. The branching property ensures that $(Y_{t}(x+y), t\geq 0)\overset{law}{=}(Y_{t}(x)+X_{t}(y), t\geq 0)$ where $(X_{t}(y), t\geq 0)$ is a CBI$(\Psi, 0)$ starting from $y$ (that is a CB-process without immigration and with branching mechanism $\Psi$) independent of $(Y_{t}(x), t\geq 0)$. The Kolmogorov's extension theorem allows one to construct a flow $(Y_{t}(x), t\geq 0, x \geq 0)$ such that for every $y\geq 0$, $(Y_{t}(x+y)-Y_{t}(x), t\geq 0)$ has the same law as $(X_{t}(y), t\geq 0)$ a CB-process started from $y$.
We denote by $(M_{t},t\geq 0)$ the Stieltjes-measure associated with the increasing process $x \in [0,1] \mapsto Y_{t}(x)$. Namely, define
\begin{align*}
&M_{t}(]x,y]):=Y_{t}(y)-Y_{t}(x),  \quad 0 \leq x \leq y \leq 1.\\
&M_{t}(\{0\}):=Y_{t}(0).
\end{align*}
The process $(Y_{t}(1), t\geq 0)$ is assumed to be conservative, therefore the process $(M_{t}, t\geq 0)$ is valued in the space $\mathcal{M}_{f}$ of finite measures on $[0,1]$. By a slight abuse of notation, we denote by $(Y_{t}, t\geq 0)$ the process $(Y_{t}(1),t\geq 0)$. The framework of measure-valued processes allows us to consider an infinitely many types model. Namely each individual has initially its own type (which lies in $[0,1]$) and transmits it to its progeny. People issued from the immigration have a \textit{distinguished} type fixed at $0$. Since the types do not evolve in time, they allow us to track the ancestors at time $0$. This model can be viewed as a superprocess without spatial motion (or without mutation in population genetics vocable).
\\
Let $\mathcal{C}$ be the class of functions on $\mathcal{M}_{f}$ of the form $$F(\eta):=G\left(\langle f_{1}, \eta \rangle, ..., \langle f_{n}, \eta \rangle\right),$$
where $\langle f, \eta \rangle:=\int_{[0,1]}f(x)\eta(dx)$, $G\in C^{2}(\mathbb{R}^{n})$ and $f_{1},...,f_{n}$ are bounded measurable functions on $[0,1]$. Section 9.3 of Li's book \cite{Li} (see Theorem 9.18 p. 218) ensures that the following operator acting on the space $\mathcal{M}_{f}$ is an extended generator of $(M_{t},t\geq 0)$. For any $\eta\in \mathcal{M}_{f}$, 
\begin{align}
\mathcal{L}F(\eta)&:=\sigma^{2}/2\int_{0}^{1}\int_{0}^{1}\eta(da)\delta_{a}(db)F''(\eta;a,b)\\
&+\beta F'(\eta;0)\\
&+\int_{0}^{1}\eta(da)\int_{0}^{\infty}\hat{\nu_{1}}(dh)[F(\eta+h\delta_{a})-F(\eta)-hF'(\eta,a)]\\
&+\int_{0}^{\infty}\hat{\nu_{0}}(dh)[F(\eta+h\delta_{0})-F(\eta)]
\end{align}
where $F'(\eta;a):=\lim_{\epsilon \rightarrow 0}\frac{1}{\epsilon}[F(\eta+\epsilon \delta_{a})-F(\eta)]$ is the Gateaux derivative of $F$ at $\eta$ in direction $\delta_{a}$, and $F''(\eta;a,b):=G'(\eta; b)$ with $G(\eta)=F'(\eta; a)$.
The terms (1) and (3) correspond to the reproduction, see for instance Section 6.1 p. 106 of Dawson \cite{Dawson0}. The terms (2) and (4) correspond to the immigration. We stress that in our model the immigration is concentrated on $0$, contrary to other works which consider infinitely many types for the immigrants.
For the interested reader, the operator $\mathcal{L}$ corresponds with that given in equation (9.25) of Section 9 of Li \cite{Li} by setting $H(d\mu)=\int_{0}^{\infty}\hat{\nu}_{0}(dh)\delta_{h\delta_{0}}(d\mu)$ and $\eta=\beta\delta_{0}.$\\
\\
For all $\eta \in \mathcal{M}_{f}$, we denote by $|\eta|$ the total mass $|\eta|:=\eta([0,1])$. If $(M_{t},t \geq 0)$ is a Markov process with the above operator for generator, the process $(|M_{t}|, t\geq 0)$ is by construction a CBI. 
This is also plain from the form of the generator $\mathcal{L}$:
let $\psi$ be a twice differentiable function on $\mathbb{R}_{+}$ and define $F:\eta \mapsto \psi(|\eta|)$, we find $\mathcal{L}F(\eta)=zG_{B}\psi(z)+G_{I}\psi(z)$ for $z=|\eta|$, where 
\begin{align}
&G_{B}\psi(z)=\frac{\sigma^{2}}{2}\psi''(z)+\int_{0}^{\infty}[\psi(z+h)-\psi(z)-h\psi'(z)]\hat{\nu_{1}}(dh)\\
&G_{I}\psi(z)=\beta \psi'(z)+\int_{0}^{\infty}[\psi(z+h)-\psi(z)]\hat{\nu_{0}}(dh).
\end{align}
\subsection{Background on $M$-generalized Fleming-Viot processes with immigration} \label{GFVI}
We denote by $\mathcal{M}_{1}$ the space of probability measures on $[0,1]$. Let $c_{0}$, $c_{1}$ be two non-negative real numbers and  $\nu_{0}$, $\nu_{1}$ be two measures on $[0,1]$ such that $\int_{0}^{1}x\nu_{0}(dx)<\infty$ and $\int_{0}^{1}x^{2}\nu_{1}(dx)<\infty$. Following the notation of \cite{coaldist}, we define the couple of finite measures $M=(\Lambda_{0}, \Lambda_{1})$ such that
\begin{center}
$\Lambda_{0}(dx)=c_{0}\delta_{0}(dx)+x\nu_{0}(dx), \ \Lambda_{1}(dx)=c_{1}\delta_{0}(dx)+x^{2}\nu_{1}(dx)$.
\end{center}
The $M$-generalized Fleming-Viot process with immigration describes a population with \textit{constant size} which evolves by resampling. Let $(\rho_{t},t\geq 0)$ be an $M$-generalized Fleming-Viot process with immigration. The evolution of this process is a superposition of a continuous evolution, and a discontinuous one. The continuous evolution can be described as follows: every couple of individuals is sampled at constant rate $c_{1}$, in which case one of the two individuals gives its type to the other: this is a reproduction event. Furthermore, any individual is picked at constant rate $c_{0}$, and its type replaced by the distinguished type $0$ (the immigrant type): this is an immigration event. The discontinuous evolution is prescribed by two independent Poisson point measures $N_{0}$ and $N_{1}$ on $\mathbb{R}_{+}\times [0,1]$ with respective intensity $dt \otimes \nu_{0}(dx)$ and $dt \otimes \nu_{1}(dx)$. More precisely, if $(t,x)$ is an atom of $N_{0}+N_{1}$ then $t$ is a jump time of the process $(\rho_{t},t\geq 0)$ and the conditional law of $\rho_{t}$ given $\rho_{t-}$ is:
\begin{itemize} 
\item $(1-x)\rho_{t-}+x\delta_{U}$, if $(t,x)$ is an atom of $N_{1}$, where $U$ is distributed according to $\rho_{t-}$
\item $(1-x)\rho_{t-}+x\delta_{0}$, if $(t,x)$ is an atom of $N_{0}$.
\end{itemize}
If $(t,x)$ is an atom of $N_{1}$, an individual is picked at random in the population at generation $t-$ and generates a proportion $x$ of the population at time $t$: this is a reproduction event, as for the genuine generalized Fleming-Viot process (see \cite{LGB1} p278).
If $(t,x)$ is an atom of $N_{0}$, the individual $0$ at time $t-$ generates a proportion $x$ of the population at time $t$: this is an immigration event. In both cases, the population at time $t-$ is reduced by a factor $1-x$ so that, at time $t$, the total size is still $1$. The genealogy of this population (which is identified as a probability measure on $[0,1]$) is given by an $M$-coalescent (see Section \ref{genealogy} below). This description is purely heuristic (we stress for instance that the atoms of $N_{0}+N_{1}$ may form an infinite dense set), to make a rigorous construction of such processes, we refer to the Section 5.2 of \cite{coaldist} (or alternatively Section 3.2 of \cite{partitionflow}). \\
For any $p\in \mathbb{N}$ and any continuous function $f$ on $[0,1]^{p}$, we denote by $G_{f}$ the map 
$$\rho\in \mathcal{M}_{1}\mapsto \langle f, \rho^{\otimes p} \rangle := \int_{[0,1]^{p}}f(x)\rho^{\otimes p} (dx) =\int_{[0,1]^{p}}f(x_{1},...,x_{p})\rho(dx_{1})...\rho(dx_{p}).$$ Let $(\mathcal{F},\mathcal{D})$ denote the generator of $(\rho_{t},t\geq 0)$ and its domain. The vector space generated by the functionals of the type $G_{f}$ forms a core of $(\mathcal{F},\mathcal{D})$ and we have (see Lemma 5.2 in \cite{coaldist}):
\begin{align*}
\mathcal{F}G_{f}(\rho)&=c_{1}\sum_{1\leq i<j\leq p}\int_{[0,1]^{p}}[f(\textsl{x}^{i,j})-f(\textsl{x})]\rho^{\otimes p}(d\textsl{x}) \tag{1'} \\
&+c_{0}\sum_{1\leq j\leq p}\int_{[0,1]^{p}}[f(\textsl{x}^{0,j})-f(\textsl{x})]\rho^{\otimes p}(d\textsl{x}) \tag{2'}\\
&+\int_{0}^{1} \nu_{1}(dr)\int \rho(da)[G_{f}((1-r)\rho+r\delta_{a})-G_{f}(\rho)] \tag{3'}\\
&+\int_{0}^{1} \nu_{0}(dr)[G_{f}((1-r)\rho+r\delta_{0})-G_{f}(\rho)] \tag{4'}.
\end{align*}
where $\textsl{x}$ denotes the vector $(x_{1},...,x_{p})$ and
\begin{itemize}
\item the vector $\textsl{x}^{0,j}$ is defined by $\textsl{x}^{0,j}_{k}=x_{k}$, for all $k\neq j$ and $\textsl{x}^{0,j}_{j}=0$,
\item the vector $\textsl{x}^{i,j}$ is defined by $\textsl{x}^{i,j}_{k}=x_{k}$, for all $k \neq j$ and $\textsl{x}^{i,j}_{j}=x_{i}.$
\end{itemize}  

\section{Relations between CBIs and $M$-GFVIs}
\label{mesureproc}
 
\subsection{Forward results}
\label{forwardresults}

The expressions of the generators of $(M_{t},t\geq 0)$ and $(\rho_{t}, t\geq 0)$ lead us to specify the connections between CBIs and GFVIs. We add a cemetery point $\Delta$ to the space $\mathcal{M}_{1}$ and define $(R_{t}, t\geq 0):=(\frac{M_{t}}{|M_{t}|}, t\geq 0)$, the ratio process with lifetime  $\tau:=\inf \{t\geq 0 ; |M_t|=0\}$. By convention, for all $t\geq \tau$, we set $R_{t}=\Delta$. 
As mentioned in the Introduction, we shall focus our study on the two following critical CBIs:
\begin{itemize}
\item[(i)] $(Y_{t}, t\geq 0)$ is CBI with parameters $\sigma^{2},\beta \geq 0$ and $\hat{\nu}_{0}=\hat{\nu}_{1}=0$, so that $\Psi(q)=\frac{\sigma^{2}}{2}q^{2}$ and $\Phi(q)=\beta q$.
\item[(ii)] $(Y_{t}, t\geq 0)$ is a CBI with $\sigma^{2}=\beta=0$, $\hat{\nu_{0}}(dh)=c'h^{-\alpha}1_{h>0}dh$ and $\hat{\nu_{1}}(dh)=ch^{-1-\alpha}1_{h>0}dh$ for $1<\alpha<2$, so that
$\Psi(q)=dq^{\alpha}$ and $\Phi(q)=d'\alpha q^{\alpha-1}$ with $d'=\frac{\Gamma(2-\alpha)}{\alpha(\alpha-1)}c'$ and $d=\frac{\Gamma(2-\alpha)}{\alpha(\alpha-1)}c$
\end{itemize}
Notice that the CBI in (i) may be seen as a limit case of the CBIs in (ii) for $\alpha= 2$.
We first establish in the following proposition a dichotomy for the finiteness of the lifetime, depending on the ratio immigration over reproduction. 
\begin{prop} \label{extCBI}
Recall the notation $\tau=\inf \{t\geq 0, Y_{t}=0\}.$
\begin{itemize} 
\item If $\frac{\beta}{\sigma^{2}}\geq \frac{1}{2}$ in case (i) or $\frac{c'}{c}\geq \frac{\alpha-1}{\alpha}$ in case (ii), then $\mathbb{P}[\tau=\infty]=1$. 
\item If $\frac{\beta}{\sigma^{2}}< \frac{1}{2}$ in case (i) or $\frac{c'}{c}< \frac{\alpha-1}{\alpha}$ in case (ii), then $\mathbb{P}[\tau<\infty]=1$.
\end{itemize}
\end{prop}
We then deal with the random change of time. In the case of a CB-process (that is a CBI process without immigration), Birkner \textit{et al.} used the Lamperti representation and worked on the embedded stable spectrally positive L\'evy process. We shall work directly on the CBI process instead. 
For $0 \leq t \leq \tau$, we define:
$$ C(t) = \int_{0}^{t}Y_{s}^{1-\alpha}ds,$$
in case (ii) and set $\alpha=2$ in case (i). 
\begin{prop} \label{timechange}
In both cases (i) and (ii), 
we have:
$$\mathbb{P}\left(C(\tau)=\infty\right)=1.$$
In other words, the additive functional $C$ maps $[0, \tau[$ to $[0, \infty[$.
\end{prop}
By convention, if $\tau$ is almost surely finite we set $C(t)=C(\tau)=\infty$ for all $t\geq \tau$. Denote by $C^{-1}$ the right continuous inverse of the functional $C$. This maps $[0,\infty[$ to $[0, \tau[$,
a.s. We stress that in most cases, $(R_{t}, t\geq 0)$ is not a Markov process. Nevertheless, in some cases, through a change of time, the process $(R_{t}, t\geq 0)$ may be changed into a Markov process. This shall be stated in the following Theorem where the functional $C$ is central.
\\
For every $x, y >0$, denote by $Beta(x,y)(dr)$ the finite measure with density $$r^{x-1}(1-r)^{y-1}1_{(0,1)}(r)dr,$$ and recall that its total mass is given by the Beta function $B(x,y)$.

\begin{thm} \label{thm1} Let $(M_{t}, t\geq 0)$ be the measure-valued process associated to a process $(Y_{t}(x), x\in [0,1], t\geq 0).$
\begin{itemize}
\item[-]  In case (i), the process $(R_{C^{-1}(t)})_{t\geq 0}$ is a $M$-Fleming-Viot process with immigration with
\begin{center}
$\Lambda_{0}(dr)=\beta\delta_{0}(dr)$ and $\Lambda_{1}(dr)=\sigma^{2}\delta_{0}(dr)$.
\end{center}
\item[-] In case (ii), the process $(R_{C^{-1}(t)})_{t\geq 0}$ is a $M$-generalized Fleming-Viot process with immigration with
\begin{center}
$\Lambda_{0}(dr)=c'Beta(2-\alpha,\alpha-1)(dr)$ and $\Lambda_{1}(dr)=cBeta(2-\alpha,\alpha)(dr)$.
\end{center} 
\end{itemize}
\end{thm}

The proof requires rather technical arguments on the generators and is given in Section \ref{Proofthm1}.
\begin{remark}
\begin{itemize} 
\item The CBIs in the statement of Theorem \ref{thm1} with $\sigma^2=\beta$ in case $(i)$ or $c=c'$ in case $(ii)$, are also CBs conditioned on non extinction and are studied further in Section \ref{genealogy}.
\item Contrary to the case without immigration, see Theorem 1.1 in \cite{Birk}, we have to restrict ourselves to $\alpha \in (1,2]$.

\end{itemize}
\end{remark}
So far, we state that the ratio process $(R_{t}, t\geq 0)$ associated to $(M_{t}, t\geq 0)$, once time changed by $C^{-1}$, is a $M$-GFVI process. Conversely, starting from a $M$-GFVI process, we could wonder how to recover the measure-valued CBI process $(M_{t}, t\geq 0)$. This lead us to investigate the relation between the time changed ratio process $(R_{C^{-1}(t)}, t\geq 0)$ and the process $(Y_{t}, t\geq 0)$. 
\begin{prop} \label{independence} In case $(i)$ of Theorem \ref{thm1}, the additive functional $(C(t), t\geq 0)$ and $(R_{C^{-1}(t)}, 0\leq t< \tau)$ are independent.
\end{prop}
This proves that in case $(i)$ we need additional randomness to reconstruct $M$ from the $M$-GFVI process. On the contrary, in case $(ii)$, the process $(Y_{t}, t\geq 0)$ is clearly not independent of the ratio process $(R_{t}, t\geq 0)$,
since both processes jump at the same time. 
\\
The proof of Propositions \ref{extCBI}, \ref{timechange} are given in the next Subsection. Some rather technical arguments are needed to prove Proposition \ref{independence}. We postpone its proof to the end of Section \ref{Proofthm1}. 
\subsection{Proofs of Propositions \ref{extCBI}, \ref{timechange}}
\label{Proof0}
\textit{Proof of Proposition \ref{extCBI}.} 
Let $(X_{t}(x),t\geq 0)$ denote an $\alpha$-stable branching process started at $x$ (with $\alpha \in (1,2]$). Denote $\zeta$ its absorption time, $\zeta:=\inf \{t\geq 0; X_{t}(x)=0\}$. The following construction of the process $(Y_{t}(0),t\geq 0)$ may be deduced from the expression of the Laplace transform of the CBI process. 
We shall need the canonical measure $\N$ which is a sigma-finite measure on c\`adl\`ag paths and represents informally the ``law'' of the population generated by one single individual in a CB($\Psi$), see Li \cite{Li}. 
We write: 
\begin{equation}
\label{decompositionpoisson}
(Y_{t}(0),t\geq 0)=\left(\sum_{i\in \mathcal{I}} X^i_{(t-t_{i})_{+}}, t \geq0 \right) 
\end{equation}
with $\sum_{i} \delta_{(t_i,X^i)}$ a Poisson random measure on $\mathbb{R}_{+}\times \mathcal{D}(\mathbb{R}_{+}, \mathbb{R}_{+})$ with intensity $dt \otimes \mu$, where $\mathcal{D}(\mathbb{R}_{+}, \mathbb{R}_{+})$ denotes the space of c\`adl\`ag functions, and $\mu$ is defined as follows:
\begin{itemize}
\item in case $(ii)$, $\mu(dX)=\int \hat{\nu_0}(dx) \mathbb{P}_{x}(dX)$, where $\mathbb{P}_x$ is the law of a CB($\Psi$) with $\Psi(q)=dq^{\alpha}$.
Formula (\ref{decompositionpoisson}) may be understood as follows: at the jump times $t_i$ of a pure jump stable subordinator with L\'evy measure $\hat{\nu_0}$, a new arrival of immigrants, of size $X^i_0$, occurs in the population.
Each of these "packs", labelled by $i\in \mathcal{I}$, generates its own descendance $(X^i_t, t \geq 0)$, which is a CB($\Psi$) process.

\item in case $(i)$, $\mu(dX)=\beta \ \N(dX)$, where $\N$ is the canonical measure associated to the CB($\Psi$) with $\Psi(q)=\frac{\sigma^{2}}{2}q^{2}$. The canonical measure may be thought of as the ``law'' of the population generated by one single individual.
The link with case (ii) is the following:
the pure jump subordinator degenerates into a continuous subordinator equal to $( t \mapsto \beta t)$. The immigrants no more arrive by packs, but appear continuously.
\end{itemize}
Actually, the canonical measure $\N$ is defined in both cases (i) and (ii), and 
we may always write $\mu(dX)= \Phi(\N(dX)).$ The process $(Y_{t}(0), t\geq 0)$ is a CBI$(\Psi, \Phi)$ started at $0$. We call $\mathcal{R}$ the set of zeros of $(Y_{t}(0), t> 0)$:  \[\mathcal{R}:=\{t > 0; Y_{t}(0)=0\}.\]
Denote $\zeta_{i}= \inf{\{ t >0, X^i_t=0 \}}$ the lifetime of the branching process $X^{i}$. The intervals $]t_{i}, t_{i}+\zeta_{i}[$ and $[t_{i}, t_{i}+\zeta_{i}[$ represent respectively the time where $X^{i}$ is alive in case (i) and in case (ii) (in this case, we have $X^{i}_{t_{i}}>0$.) 
Therefore, if we define $\mathcal{\tilde{R}}$ as the set of the positive real numbers left uncovered by the random intervals 
$]t_{i},t_{i}+\zeta_{i}[$, that is:
\[\mathcal{\tilde{R}}:= \mathbb{R}^{\star}_{+}\setminus \bigcup_{i\in \mathcal{I}} \ ]t_{i},t_{i}+\zeta_{i}[. \] 
we have $\mathcal{R} \subset \mathcal{\tilde{R}}$ with equality in case (i) only.

%

The lengths $\zeta_i$ have law $\mu(\zeta\in dt)$ thanks to the Poisson construction of $Y(0)$. We now distinguish the two cases:  
\begin{itemize}
\item Feller case: this corresponds to $\alpha=2$. We have $\Psi(q):=\frac{\sigma^{2}}{2}q$ and $\Phi(q):=\beta q$, and thus
\[\mu[\zeta>t]=\beta \; \mathbb{N}[\zeta>t]=\frac{2 \beta}{\sigma^{2}}\frac{1}{t}\] see Li \cite{Li} p. 62. Using Example 1 p. 180 of Fitzsimmons et al. \cite{Fitzsimmons}, we deduce that
\begin{eqnarray} \label{ineq1}
\mathcal{\tilde{R}}=\emptyset \   \text{ a.s. if and only if  }\  \frac{2 \beta}{\sigma^2}\geq 1.
\end{eqnarray}
\item Stable case: this corresponds to $\alpha \in  (1,2)$. Recall $\Psi(q):=dq^{\alpha}, \Phi(q):=d'\alpha q^{\alpha-1}$. In that case, we have, 
\[\N(\zeta>t)=d^{-\frac{1}{\alpha-1}}[(\alpha-1)t]^{-\frac{1}{\alpha-1}}. \]
 Thus, $ \mu[\zeta>t]=\Phi(\N(\zeta>t))= \frac{\alpha}{\alpha-1} \frac{d'}{d}\frac{1}{t}$. Recall that $\frac{d'}{d}=\frac{c'}{c}$. Therefore, using reference \cite{Fitzsimmons}, we deduce that
\begin{eqnarray} \label{ineq2}
\mathcal{\tilde{R}}=\emptyset \  \text{ a.s. if and only if } \ \frac{c'}{c}\geq \frac{\alpha-1}{\alpha}.
\end{eqnarray}
\end{itemize}
This allows us to establish the first point of Proposition \ref{extCBI}: 
we get $\mathcal{R} \subset \mathcal{\tilde{R}}= \emptyset$, and the inequality $Y_{t}(1)\geq Y_{t}(0)$ for all $t$ ensures that $\tau=\infty$.\\ 

We deal now with the second point of Proposition \ref{extCBI}. Assume that $\frac{c'}{c}<\frac{\alpha-1}{\alpha}$ or $\frac{\beta}{\sigma^{2}}<\frac{1}{2}$. By assertions (\ref{ineq1}) and (\ref{ineq2}), we already know that $\mathcal{\tilde{R}}\neq \emptyset$. 
However, what we really need is that $\mathcal{\tilde{R}}$ is a.s. not bounded.
To that aim, observe that, 
in both cases (i) and (ii),
 \[\mu[\zeta>s]=\Phi(\N(\zeta >s))= \frac{\kappa}{s}\] with $\kappa=\frac{\alpha}{\alpha-1} \frac{d'}{d}=\frac{\alpha}{\alpha-1} \frac{c'}{c}<1$
if $1 < \alpha <2$ and $\kappa= \frac{2 \beta}{\sigma^2}<1$ if $\alpha=2$.
Thus $\int_{1}^{u}\mu[\zeta>s]ds=\kappa \ln(u)$ and we obtain \[\exp\left( -\int_{1}^{u}\mu[\zeta>s]ds \right)=\left(\frac{1}{u}\right)^{\kappa}.\] 
Therefore, since $\kappa<1$, $$\int_{1}^{\infty}\exp\left( -\int_{1}^{u}\mu[\zeta>s]ds \right)du= \infty,$$ 
which implies thanks to Corollary 4 (Equation 17 p 183) of \cite{Fitzsimmons} that $ \mathcal{\tilde{R}}$ is a.s. not bounded. 
\\

Since $\mathcal{R}=\tilde{\mathcal{R}}$ in case (i), the set $\mathcal{R}$ is a.s. not bounded in that case. Now, we prove that $\mathcal{R}$ is a.s. not bounded in case (ii). The set $\mathcal{\tilde{R}}$ is almost surely not empty and not bounded. Moreover this is a perfect set (Corollary 1 of \cite{Fitzsimmons}). Since there are only countable points $(t_{i}, i\in \mathcal{I})$, the set $\tilde{\mathcal{R}}=\mathcal{R}\setminus \bigcup_{i\in I}\{t_{i}\}$ is also uncountable and not bounded.
\\	

Last, recall from Subsection \ref{CBI} that we may write $Y_{t}(1)= Y_{t}(0)+X_{t}(1)$ for all $t\geq 0$ with $(X_{t}(1), t\geq 0)$ a CB-process independent of $(Y_{t}(0), t\geq 0)$. Let $\xi:=\inf \{t\geq 0, X_{t}(1)=0\}$ be the extinction time of $(X_{t}(1), t\geq 0)$. Since $\mathcal{R}$ is a.s. not bounded in both cases (i) and (ii), $\mathcal{R} \cap (\xi,\infty) \neq \emptyset$, and  $\tau <\infty$ almost surely. 
$\square$
\\
\\
\textit{Proof of Proposition \ref{timechange}.} 
Recall that $Y_{t}(x)$ is the value of the CBI started at $x$ at time $t$. We will denote by $\tau^{x}(0):=\inf{\{t>0, Y_{t}(x)=0\}}$. 
With this notation, $\tau^{1}(0)=\tau$ introduced in Section \ref{forwardresults}.
In both cases (i) and (ii), the processes are self-similar, see Kyprianou and Pardo \cite{Kyprianou}. Namely, we have $$\left(x Y_{x^{1-\alpha}t}(1), t \geq 0\right) \overset{law}{=} \left(Y_{t}(x), t \geq 0 \right),$$ where we take $\alpha=2$ in case $(i)$. Performing the change of variable $s=x^{1-\alpha} t$, we obtain
\begin{align}
\label{equalityinlaw}
\int_{0}^{\tau^{x}(0)} dt \ Y_{t}(x)^{1-\alpha}  &\overset{law}{=}  \int_{0}^{\tau^{1}(0)} ds \  Y_{s}(1)^{1-\alpha}.  
\end{align}
According to Proposition \ref{extCBI}, depending on the values of the parameters:
\begin{itemize}
 \item Either $\P(\tau^{x}(0)<\infty)=1$ for every $x$. Let $x>1$. Denote $\tau^{x}(1)=\inf{ \{ t>0, Y_{t}(x) \leq 1 \}}$. We have $\P(\tau^{x}(1)<\infty)=1$. We have: 
$$\int_{0}^{\tau^{x}(0)} dt \ Y_{t}(x)^{1-\alpha} = \int_{0}^{\tau^{x}(1)} dt \ Y_{t}(x)^{1-\alpha} + \int_{\tau^{x}(1)}^{\tau^{x}(0)} dt \ Y_{t}(x)^{1-\alpha}$$ 
By the strong Markov property applied at the stopping time $\tau^{x}(1)$, since $Y$ has no negative jumps:
\[\int_{\tau^{x}(1)}^{\tau^{x}(0)} dt \ Y_{t}(x)^{1-\alpha}\overset{law}{=}\int_{0}^{\tau^{1}(0)} dt \ \tilde{Y}_{t}(1)^{1-\alpha},\]
with $(\tilde{Y}_{t}(1),t\geq 0)$ an independent copy started from $1$. 
Since $$\int_{0}^{\tau^{x}(1)} dt \ Y_{t}(x)^{1-\alpha} >0, \ a.s. ,$$ the equality \reff{equalityinlaw} is impossible unless both sides of the equality are infinite almost surely. We thus get that $C(\tau)=\infty$ almost surely in that case. 
 \item Either $\P(\tau^{x}(0)=\infty)=1$ for every $x$, on which case we may rewrite \reff{equalityinlaw} as follows:
$$ \int_{0}^{\infty} dt \ Y_{t}(x)^{1-\alpha}  \overset{law}{=}  \int_{0}^{\infty} ds \  Y_{s}(1)^{1-\alpha}.$$
Since, for $x >1$, the difference $(Y_{t}(x)-Y_{t}(1), t \geq 0)$ is an $\alpha$-stable CB-process started at $x-1>0$, we deduce that $C(\tau)=\infty$ almost surely again. 
\end{itemize}
This proves the statement. $\square$
\begin{remark} \label{remark}
The situation is quite different when the CBI process starts at $0$, in which case the time change also diverges in the neighbourhood of $0$. The same change of variables as in (\ref{equalityinlaw}) yields, for all $0<x<k$,
\begin{align*}
\int_{0}^{\iota^{x}(k)} dt \ Y_{t}(x)^{1-\alpha}  &\overset{law}{=}  \int_{0}^{\iota^1(k/x)} dt \  Y_{t}(1)^{1-\alpha},
\end{align*}
with $\iota^{x}(k)= \inf\{t>0, Y_{t}(x) \geq k \} \in [0,\infty]$. Letting $x$ tend to $0$, we get $\iota^{1}(k/x)\longrightarrow \infty$ and the right hand side diverges to infinity. Thus, the left hand side also diverges, which implies that: 
$$\P\left(\int_{0}^{\iota^{0}(k)} dt \ Y_{t}(0)^{1-\alpha} = \infty\right)=1.$$
\end{remark}
\section{Genealogy of the Beta-Fleming-Viot processes with immigration}  \label{genealogy}
To describe the genealogy associated with stable CBs, Bertoin and Le Gall \cite{LGB3} and Birkner et al. \cite{Birk} used partition-valued processes called Beta-coalescents. These processes form a subclass of $\Lambda$-coalescents, introduced independently by Pitman and Sagitov in 1999. A $\Lambda$-coalescent is an exchangeable process in the sense that its law is invariant under the action of any permutation. In words, there is no distinction between the individuals. Although these processes arise as models of genealogy for a wide range of stochastic populations, they are not in general adapted to describe the genealogy of a population with immigration. Recently, a larger class of processes called $M$-coalescents has been defined in \cite{coaldist} (see Section 5). These processes are precisely those describing the genealogy of $M$-GFVIs.
\begin{remark}
We mention that the use of the lookdown construction in Birkner et al. \cite{Birk} may be easily adapted to our framework and yields a genealogy for any conservative CBI. Moreover, other genealogies, based on continuous trees, have been investigated by Lambert \cite{Lambert} and Duquesne \cite{Duquesne2}.
\end{remark}
\subsection{Background on $M$-coalescents}
Before focusing on the $M$-coalescents involved in the context of Theorem \ref{thm1}, we recall their general definition and the duality with the $M$-GFVIs. Contrary to the $\Lambda$-coalescents, the $M$-coalescents are only invariant by permutations 
letting $0$ fixed.
The individual $0$ represents the immigrant lineage and is distinguished from the others. We denote by $\mathcal{P}^{0}_{\infty}$ the space of partitions of $\mathbb{Z}_{+}:=\{0\}\bigcup \mathbb{N}$. Let $\pi \in \mathcal{P}^{0}_{\infty}$. By convention, we identify $\pi$ with the sequence $(\pi_{0},\pi_{1},...)$ of the blocks of $\pi$ enumerated in increasing order of their smallest element: for every $i\leq j$, $\min \pi_{i}\leq \min \pi_{j}$. Let $\cro{n}$ denote the set $\{0,...,n\}$ and $\mathcal{P}^{0}_{n}$ the space of partitions of $\cro{n}$. The partition of $\cro{n}$ into singletons is denoted by $0_{\cro{n}}$. As in Section 2.2, the notation $M$ stands for a pair of finite measures $(\Lambda_{0}, \Lambda_{1})$ such that:
\begin{center}
$\Lambda_{0}(dx)=c_{0}\delta_{0}(dx)+x\nu_{0}(dx), \ \Lambda_{1}(dx)=c_{1}\delta_{0}(dx)+x^{2}\nu_{1}(dx)$,
\end{center}
where $c_{0}$, $c_{1}$ are two non-negative real numbers and $\nu_{0}$, $\nu_{1}$ are two measures on $[0,1]$ subject to the same conditions as in Section \ref{GFVI}. Let $N_{0}$ and $N_{1}$ be two Poisson point measures with intensity respectively $dt\otimes \nu_{0}$ and $dt\otimes \nu_{1}$. An $M$-coalescent is a Feller process $(\Pi(t), t\geq 0)$ valued in $\mathcal{P}^{0}_{\infty}$ with the following dynamics.
\begin{itemize}
\item At an atom $(t,x)$ of $N_{1}$, flip a coin with probability of "heads" $x$ for each block not containing $0$. All blocks flipping "heads" are merged immediately in one block. At time $t$, a proportion $x$ share a common parent in the population. 
\item At an atom $(t,x)$ of $N_{0}$, flip a coin with probability of "heads" $x$ for each block not containing $0$. All blocks flipping "heads" coagulate immediately with the distinguished block. At time $t$, a proportion $x$ of the population is children of immigrant.
\end{itemize}
In order to take into account the parameters $c_{0}$ and $c_{1}$, imagine that at constant rate $c_{1}$, two blocks (not containing $0$) merge \textit{continuously} in time, and at constant rate $c_{0}$, one block (not containing $0$) merged with the distinguished one. We refer to Section 4.2 of \cite{coaldist} for a rigorous definition. Let $\pi \in \mathcal{P}^{0}_{n}$. The jump rate of an $M$-coalescent from $0_{\cro{n}}$ to $\pi$, denoted by $q_{\pi}$, is given as follows:
\begin{itemize}
\item If $\pi$ has one block not containing $0$ with $k$ elements and $2\leq k\leq n$, then $$q_{\pi}=\lambda_{n,k}:=\int_{0}^{1}x^{k-2}(1-x)^{n-k}\Lambda_{1}(dx).$$
\item If the distinguished block of $\pi$ has $k+1$ elements (counting $0$) and $1\leq k\leq n$ then $$q_{\pi}=r_{n,k}:=\int_{0}^{1}x^{k-1}(1-x)^{n-k}\Lambda_{0}(dx).$$
\end{itemize}The next duality property is a key result and links the $M$-GFVIs to the $M$-coalescents. For any $\pi$ in $\mathcal{P}^{0}_{\infty}$, define
\begin{center}
$\alpha_{\pi} : k\mapsto$ the index of the block of $\pi$ containing $k$. \end{center}
We have the duality relation (see Lemma 4 in \cite{partitionflow}): for any $p\geq 1$ and $f\in C([0,1]^{p})$,
$$\mathbb{E}\left[\int_{[0,1]^{p+1}}f(x_{\alpha_{\Pi(t)}(1)},...,x_{\alpha_{\Pi(t)}(p)})\delta_{0}(dx_{0}) dx_{1}...dx_{p}\right]=\mathbb{E}\left[\int_{[0,1]^{p}}f(x_{1},...,x_{p})\rho_{t}(dx_{1})...\rho_{t}(dx_{p})\right],$$
where $(\rho_{t}, t\geq 0)$ is a $M$-GFVI started from the Lebesgue measure on $[0,1]$. 
We establish a useful lemma relating genuine $\Lambda$-coalescents and $M$-coalescents. Consider a $\Lambda$-coalescent taking values in the set $\mathcal{P}^{0}_{\infty}$; this differs from the usual convention, according to which they are valued in the set $\mathcal{P}_{\infty}$ of the partitions of $\N$ (see Chapters 1 and 3 of \cite{Beres2} for a complete introduction to these processes). In that framework, $\Lambda$-coalescents appear as a subclass of $M$-coalescents and the integer $0$ may be viewed as a typical individual. 
The proof is postponed in Section \ref{proofs}.
\begin{lem} \label{Lambda}
A $M$-coalescent, with $M=(\Lambda_{0},\Lambda_{1})$ is also a $\Lambda$-coalescent on $\mathcal{P}^{0}_{\infty}$ if and only if 
$$(1-x)\Lambda_{0}(dx)=\Lambda_{1}(dx).$$
In that case $\Lambda=\Lambda_{0}$.
\end{lem}
\subsection{The $Beta(2-\alpha, \alpha-1)$-coalescent}
The aim of this Section is to show how a $Beta(2-\alpha, \alpha-1)$-coalescent is embedded in the genealogy of an $\alpha$-stable CB-process conditioned to be never extinct. Along the way, we also derive the fixed time genealogy of the Feller CBI.\\

We first state the following straightforward Corollary of Theorem \ref{thm1}, which gives the genealogy of the ratio process at the random time $C^{-1}(t)$:
\begin{cor} 
\label{cordual}
Let $(R_{t},t\geq 0)$ be the ratio process of a CBI in case $(i)$ or $(ii)$. 
We have for all $t\geq 0$:
\begin{equation*}
\mathbb{E}\left[\int_{[0,1]^{p+1}}f(x_{\alpha_{\Pi(t)}(1)},...,x_{\alpha_{\Pi(t)}(p)})\delta_{0}(dx_{0}) dx_{1}...dx_{p}\right]=\mathbb{E}\left[\int_{[0,1]^{p}}f(x_{1},...,x_{p})R_{C^{-1}(t)}(dx_{1})...R_{C^{-1}(t)}(dx_{p})\right],
\end{equation*}
where:
\begin{itemize}
\item In case $(i)$, $(\Pi(t), t\geq 0)$ is a $M$-coalescent with $M=(\beta\delta_{0}, \sigma^{2}\delta_{0})$,
\item In case $(ii)$, $(\Pi(t), t\geq 0)$ is a $M$-coalescent with $M=(c'Beta(2-\alpha,\alpha-1),cBeta(2-\alpha,\alpha))$.  
\end{itemize}
\end{cor}

In general, we cannot set the random quantity $C(t)$ instead of $t$ in the equation of Corollary \ref{cordual}. Nevertheless, using the independence property proved in Proposition \ref{independence}, we get the following Corollary, whose proof may be found in Section \ref{proofs}..
\begin{cor} 
\label{cordual2}
In case (i),
assume $\frac{\beta}{\sigma^2} \geq \frac{1}{2}$, then for all $t\geq 0$,
$$\mathbb{E}\left[\int_{[0,1]^{p+1}}f(x_{\alpha_{\Pi(C(t))}(1)},...,x_{\alpha_{\Pi(C(t))}(p)})\delta_{0}(dx_{0}) dx_{1}...dx_{p}\right]=\mathbb{E}\left[\int_{[0,1]^{p}}f(x_{1},...,x_{p})R_{t}(dx_{1})...R_{t}(dx_{p})\right],$$
where $(\Pi(t), t\geq 0)$ is a $M$-coalescent with $M=(\beta\delta_{0}, \sigma^{2}\delta_{0})$, $(Y_{t}, t\geq 0)$ is a CBI in case (i) independent of $(\Pi(t), t\geq 0)$ and $(C(t), t\geq 0)=\left( \int_{0}^{t} \frac{1}{Y_s} ds, t\geq 0\right)$.
\end{cor}
We stress on a fundamental difference between Corollaries \ref{cordual} and \ref{cordual2}. Whereas the first gives the genealogy of the ratio process $R$ at the random time $C^{-1}(t)$, the second gives the genealogy of the ratio process $R$ at a fixed time $t$. Notice that we impose the additional assumption that $\frac{\beta}{\sigma^2} \geq \frac{1}{2}$ in Corollary \ref{cordual2} for ensuring that the lifetime is infinite. Therefore, $R_t \neq \Delta$ for all $t\geq 0$, and we may consider its genealogy. 

We easily check that the $M$-coalescents for which $M=(\sigma^2 \delta_{0},\sigma^2 \delta_{0})$ and $M=(c Beta(2-\alpha,\alpha-1), c Beta(2-\alpha,\alpha))$ fulfill the conditions of Lemma \ref{Lambda}. 
Recall from Section \ref{forwardresults} the definitions of the CBIs in case (i) and (ii) .
\begin{thm} \label{coalescent} \vspace*{1mm}
\begin{itemize}
\item[(i)] If the process $(Y_{t}, t\geq 0)$ is a CBI such that $\sigma^{2}=\beta>0$, $\hat{\nu}_{1}=\hat{\nu}_{0}=0$, then the process $(\Pi(t/\sigma^2), t\geq 0)$ defined in Corollary \ref{cordual} is a Kingman's coalescent valued in $\mathcal{P}^{0}_{\infty}$.
\item[(ii)] If the process $(Y_{t}, t\geq 0)$ is a CBI such that $\sigma^{2}=\beta=0$ and $\hat{\nu}_{0}(dh)=ch^{-\alpha}dh$, $\hat{\nu}_{1}(dh)=ch^{-\alpha-1}dh$ for some constant $c>0$ then the process $(\Pi(t/c), t\geq 0)$ defined in Corollary \ref{cordual} is a $Beta(2-\alpha, \alpha-1)$-coalescent valued in $\mathcal{P}^{0}_{\infty}$.
\end{itemize}
\end{thm}
In both cases, the process $(Y_{t}, t\geq 0)$ involved in that Theorem may be interpreted as a CB-process $(X_{t}, t\geq 0)$ without immigration ($\beta=0$ or $c'=0$) conditioned on non-extinction, see Lambert \cite{Lambert2}.
We then notice that both the genealogies of the time changed Feller diffusion and of the time changed Feller diffusion conditioned on non extinction are given by the same Kingman's coalescent. On the contrary, the genealogy of the time changed $\alpha$-stable CB-process is a $Beta(2-\alpha, \alpha)$-coalescent, whereas the genealogy of the time changed $\alpha$-stable CB-process conditioned on non-extinction is a $Beta(2-\alpha, \alpha-1)$-coalescent. 
We stress that for any $\alpha \in (1,2)$ and any borelian $B$ of $[0,1]$, we have $Beta(2-\alpha, \alpha-1)(B)\geq Beta(2-\alpha, \alpha)(B)$. This may be interpreted as the additional reproduction events needed for the process to be never extinct.
\subsection{Proofs.}
\label{proofs}

\textit{Proof of Lemma \ref{Lambda}.} Let $(\Pi'(t),t\geq 0)$ be a $\Lambda$-coalescent on $\mathcal{P}^{0}_{\infty}$. Let $n\geq 1$, we may express the jump rate of $(\Pi'_{|\cro{n}}(t),t\geq 0)$ from $0_{\cro{n}}$ to $\pi$ by \\
\\
$q'_{\pi}=$
$\left\{ \begin{array}{l}
0 \mbox{ if } \pi \mbox{ has more than one non-trivial block}\\
\\
\int_{[0,1]}x^{k}(1-x)^{n+1-k}x^{-2}\Lambda(dx) \mbox{ if the non trivial block has } k \mbox{ elements}.\\
\end{array}\right.\\$
\\
Consider now a $M$-coalescent, denoting by $q_{\pi}$ the jump rate from $0_{\cro{n}}$ to $\pi$, we have \\
\\
$q_{\pi}=$
$\left\{ \begin{array}{l}
0 \mbox{ if } \pi \mbox{ has more than one non-trivial block}\\
\\
\int_{[0,1]}x^{k}(1-x)^{n-k}x^{-2}\Lambda_1(dx) \mbox{ if } \pi_{0}=\{0\} \mbox{ and the non trivial block has } k \mbox{ elements}\\
\\
\int_{[0,1]}x^{k-1}(1-x)^{n+1-k}x^{-1}\Lambda_{0}(dx) \mbox{ if } \#\pi_{0}=k.\\
\end{array}\right.\\$
\\
Since the law of a $\Lambda$-coalescent is entirely described by the family of the jump rates of its restriction on $\cro{n}$ from $0_{\cro{n}}$ to $\pi$ for $\pi$ belonging to $\mathcal{P}^{0}_{n}$ (see Section 4.2 of \cite{coursbertoin}), the processes $\Pi$ and $\Pi'$ have the same law if and only if for all $n\geq 0$ and $\pi\in \mathcal{P}^{0}_{n}$, we have $q_{\pi}=q'_{\pi}$, that is if and only if $(1-x)\Lambda_{0}(dx)=\Lambda_{1}(dx)$.
$\square$
\\
\\
\textit{Proof of Corollary \ref{cordual2}.} Since $C^{-1}(C(t))=t$,
\begin{align*}
\mathbb{E}\left[\int_{[0,1]^{p}}f(x_{1},...,x_{p})R_{t}(dx_{1})...R_{t}(dx_{p})\right] 
 &= \mathbb{E}\left[\int_{[0,1]^{p}}f(x_{1},...,x_{p})R_{C^{-1}(C(t))}(dx_{1})...R_{C^{-1}(C(t))}(dx_{p})\right]. 
\end{align*}
Then, using the independence between $R_{C^{-1}}$ and $C$, the right hand side above is also equal to:
\begin{align*}
 \int \P(C(t) \in ds) \ \mathbb{E}\left[\int_{[0,1]^{p}}f(x_{1},...,x_{p})R_{C^{-1}(s)}(dx_{1})...R_{C^{-1}(s)}(dx_{p})\right]. \\
\end{align*}
Using Corollary \ref{cordual} and choosing $(\Pi(t), t\geq 0)$ independent of $(C(t), t\geq 0)$, we find:
\begin{align*}
& \int \P(C(t) \in ds) \ \mathbb{E}\left[\int_{[0,1]^{p}}f(x_{1},...,x_{p})R_{C^{-1}(s)}(dx_{1})...R_{C^{-1}(s)}(dx_{p})\right] \\
&= \int \P(C(t) \in ds) \ \mathbb{E}\left[\int_{[0,1]^{p+1}}f(x_{\alpha_{\Pi(s)}(1)},...,x_{\alpha_{\Pi(s)}(p)})\delta_{0}(dx_{0}) dx_{1}...dx_{p}\right]\\
&= \mathbb{E}\left[\int_{[0,1]^{p+1}}f(x_{\alpha_{\Pi(C(t))}(1)},...,x_{\alpha_{\Pi(C(t))}(p)})\delta_{0}(dx_{0}) dx_{1}...dx_{p}\right].
\end{align*} 
$\square$
\begin{remark}
Notice the crucial r\^ole of the independence in order to establish Corollary \ref{cordual2}. When this property fails, as in the 
case (ii), 
the question of describing the fixed time genealogy of the $\alpha$-stable CB or CBI remains open. We refer to the discussion in Section 2.2 of Berestycki \textit{et. al} \cite{Beres3}. 
\end{remark}
\section{Proof of Theorem \ref{thm1} and Proposition \ref{independence}} \label{Proofthm1}
We first deal with Theorem \ref{thm1}. The proof of Proposition \ref{independence} is rather technical ans is postponed at the end of this Section. In order to get the connection between the two measure-valued processes $(R_{t}, t\geq 0)$ and $(M_{t}, t\geq 0)$, we may follow the ideas of Birkner \textit{et al.} \cite{Birk} and rewrite the generator of the process $(M_{t},t\geq 0)$ using the "polar coordinates": for any $\eta \in \mathcal{M}_{f}$, we define $$z:=|\eta| \mbox{ and } \rho:=\frac{\eta}{|\eta|}.$$ 
The proof relies on five lemmas. Lemma \ref{testfunction} establishes that the law of a generalized Fleming-Viot process with immigration is entirely determined by the generator $\mathcal{F}$ on the test functions of the form $\rho \mapsto \langle \phi, \rho \rangle^{m}$ with $\phi$ a measurable non-negative bounded map and $m\in \mathbb{N}$. Lemmas \ref{continuous}, \ref{genratio} and \ref{stable} allow us to study the generator $\mathcal{L}$ on the class of functions of the type $F:\eta \mapsto \frac{1}{|\eta|^{m}}\langle \phi, \eta\rangle ^{m}$. Lemma \ref{stablemes} (lifted from Lemma 3.5 of \cite{Birk}) relates stable L\'evy-measures and Beta-measures. 
We end the proof using results on time change by the inverse of an additive functional. We conclude thanks to a result due to Volkonski{\u\i} in \cite{timechange} about the generator of a time-changed process.
\begin{lem} \label{testfunction} The following martingale problem is well-posed: for any function $f$ of the form:
$$(x_{1},...,x_{p})\mapsto \prod_{i=1}^{p}\phi(x_{i})$$
with $\phi$ a non-negative measurable bounded map and $p\geq 1$, the process
$$G_{f}(\rho_{t})-\int_{0}^{t}\mathcal{F}G_{f}(\rho_{s})ds$$
is a martingale.
\end{lem}
\textit{Proof.}
Only the uniqueness has to be checked. We shall establish that the martingale problem of the statement is equivalent to the following martingale problem: for any continuous function $f$ on $[0,1]^{p}$, the process 
$$G_{f}(\rho_{t})-\int_{0}^{t}\mathcal{F}G_{f}(\rho_{s})ds$$
is a martingale. This martingale problem is well posed, see Proposition 5.2 of \cite{coaldist}. Notice that we can focus on continuous and symmetric functions since for any continuous $f$, $G_{f}=G_{\tilde{f}}$ with $\tilde{f}$ the symmetrized version of $f$. Moreover, by the Stone-Weierstrass theorem, any symmetric continuous function $f$ from $[0,1]^{p}$ to $\mathbb{R}$ can be uniformly approximated by linear combination of functions of the form $(x_{1},...,x_{p})\mapsto \prod_{i=1}^{p}\phi(x_{i})$ for some function $\phi$ continuous on $[0,1]$. We now take $f$ symmetric and continuous, and let $f_{k}$ be an approximating sequence. Plainly, we have 
$$|G_{f_{k}}(\rho)-G_{f}(\rho)|\leq ||f_{k}-f||_{\infty}$$
Assume that $(\rho_{t}, t\geq 0)$ is a solution of the martingale problem stated in the lemma. Since the map $h\mapsto G_{h}$ is linear, the process $$G_{f_{k}}(\rho_{t})-\int_{0}^{t}\mathcal{F}G_{f_{k}}(\rho_{s})ds$$ is a martingale for each $k\geq 1$. We want to prove that the process $$G_{f}(\rho_{t})-\int_{0}^{t}\mathcal{F}G_{f}(\rho_{s})ds$$ is a martingale, knowing it holds for each $f_{k}$. We will show the following convergence \begin{center} $\mathcal{F}G_{f_{k}}(\rho) \underset{k\rightarrow \infty}{\longrightarrow} \mathcal{F}G_{f}(\rho)$ uniformly in $\rho$.\end{center}
Recall expressions (1') and (2') in Subsection \ref{GFVI}, one can check that the following limits are uniform in the variable $\rho$
\[\sum_{1\leq i<j\leq p}\int_{[0,1]^{p}}[f_{k}(\textsl{x}^{i,j})-f_{k}(\textsl{x})]\rho^{\otimes p}(d\textsl{x}) \underset{k\rightarrow \infty}{\longrightarrow} \sum_{1\leq i<j\leq p}\int_{[0,1]^{p}}[f(\textsl{x}^{i,j})-f(\textsl{x})]\rho^{\otimes p}(d\textsl{x})\]
and
\[\sum_{1\leq i\leq m}\int_{[0,1]^{p}}[f_{k}(\textsl{x}^{0,i})-f_{k}(\textsl{x})]\rho^{\otimes p}(d\textsl{x}) \underset{k\rightarrow \infty}{\longrightarrow} \sum_{1\leq i\leq p}\int_{[0,1]^{p}}[f(\textsl{x}^{0,i})-f(\textsl{x})]\rho^{\otimes p}(d\textsl{x}).\]
We have now to deal with the terms (3') and (4'). In order to get that the quantity
\[\int_{0}^{1}\nu(dr)\int_{0}^{1}[G_{f_{k}}((1-r)\rho+r\delta_{a})-G_{f_{k}}(\rho)]\rho(da)\]
converges toward 
\[\int_{0}^{1}\nu(dr)\int_{0}^{1}[G_{f}((1-r)\rho+r\delta_{a})-G_{f}(\rho)]\rho(da),\]
we compute $$\langle f_{k}-f, \left((1-r)\rho+r\delta_{a} \right)^{\otimes p}\rangle -\langle f_{k}-f, \rho^{\otimes p}\rangle.$$
Since the function $f_{k}-f$ is symmetric, we may expand the $p$-fold product $\langle f_{k}-f, \left((1-r)\rho+r\delta_{a} \right)^{\otimes p}\rangle$, this yields
\begin{align*}
\langle f_{k}-f, \left((1-r)\rho+r\delta_{a} \right)^{\otimes p}&\rangle -\langle f_{k}-f, \rho^{\otimes p}\rangle\\ &=\sum_{i=0}^{p}\binom{p}{i}r^{i}(1-r)^{p-i} \left( \langle f_{k}-f,\rho^{\otimes p-i}\otimes\delta_{a}^{\otimes i} \rangle- \langle f_{k}-f,\rho^{\otimes p} \rangle\right)\\
&=pr(1-r)^{p-1}\left(\langle f_{k}-f, \rho^{\otimes p-1}\otimes\delta_{a}\rangle-\langle f_{k}-f, \rho^{\otimes p}\rangle\right)\\
&\quad + \sum_{i=2}^{p}\binom{p}{i}r^{i}(1-r)^{p-i}\left(\langle f_{k}-f,\rho^{\otimes p-i}\otimes\delta_{a}^{\otimes i} \rangle-\langle f_{k}-f,\rho^{\otimes p} \rangle\right).
\end{align*}
We use here the notation \[\langle g, \mu^{\otimes m-i}\otimes\delta_{a}^{\otimes i} \rangle := \int g(x_{1},...,x_{m-i},\underbrace{a,...,a}_{i \text{ terms}})\mu(dx_{1})...\mu(dx_{m-i}).\] 
Therefore, integrating with respect to $\rho$, the first term in the last equality vanishes and we get
\[\left\lvert \int_{0}^{1}\rho(da)\left(G_{f-f_{k}}((1-r)\rho+r\delta_{a})-G_{f-f_{k}}(\rho)\right) \right\rvert \leq 2^{p+1}||f-f_{k}||_{\infty}r^{2}\]
where $||f_{k}-f||_{\infty}$ denotes the supremum of the function $|f_{k}-f|$.
Recall that the measure $\nu_{1}$ verifies $\int_{0}^{1}r^{2}\nu_{1}(dr)<\infty$, moreover the quantity $||f_{k}-f||_{\infty}$ is bounded. Thus appealing to the Lebesgue Theorem, we get the sought-after convergence. Same arguments hold for the immigration part (4') of the operator $\mathcal{F}$. Namely we have $$|G_{f-f_{k}}((1-r)\rho+r\delta_{0})-G_{f-f_{k}}(\rho)|\leq 2^{p+1}r||f_{k}-f||_{\infty}$$
and the measure $\nu_{0}$ satisfies $\int_{0}^{1}r\nu_{0}(dr)<\infty$.
Combining our results, we obtain
$$|\mathcal{F}G_{f_{k}}(\rho)-\mathcal{F}G_{f}(\rho)|\leq C||f-f_{k}||_{\infty}$$
for a positive constant $C$ independent of $\rho$. Therefore the sequence of martingales $G_{f_{k}}(\rho_{t})-\int_{0}^{t}\mathcal{F}G_{f_{k}}(\rho_{s})ds$ converges toward $$G_{f}(\rho_{t})-\int_{0}^{t}\mathcal{F}G_{f}(\rho_{s})ds,$$ which is then a martingale. $\square$
\begin{lem} \label{continuous}
Assume that $\hat{\nu_0}=\hat{\nu_1}=0$ the generator $\mathcal{L}$ of $(M_t, t \geq 0)$ is reduced to the expressions (1) and (2):
$$\mathcal{L}F(\eta)=\sigma^{2}/2\int_{0}^{1}\int_{0}^{1}\eta(da)\delta_{a}(db)F''(\eta;a,b)+\beta F'(\eta;0)$$ 
Let $\phi$ be a measurable bounded function on $[0,1]$ and $F$ be the map $\eta \mapsto G_{f}(\rho):=\langle f,\rho^{\otimes m}\rangle$ with $f(x_{1},...,x_{p})=\prod_{i=1}^{p}\phi(x_{i})$. We have the following identity
$$|\eta| \mathcal{L}F(\eta)=\mathcal{F}G_{f}(\rho),$$
for $\eta \neq 0$, where $\mathcal{F}$ is the generator of a Fleming-Viot process with immigration with reproduction rate $c_{1}=\sigma^{2}$ and immigration rate $c_{0}=\beta$, see expressions (1') and (2').
\end{lem}
\textit{Proof.} By the calculations in Section 4.3 of Etheridge \cite{Etheridge} (but in a non-spatial setting, see also the proof of Theorem 2.1 p. 249 of Shiga \cite{Shiga}), we get:
\begin{align*}
\frac{\sigma^{2}} {2}\int_{0}^{1}\int_{0}^{1}\eta(da)\delta_{a}(db)F''(\eta;a,b)&=|\eta|^{-1}\frac{\sigma^{2}}{2}\int_{0}^{1}\int_{0}^{1}\frac{\partial^{2}G_{f}}{\partial\rho(a)\partial\rho(b)}(\rho)[\delta_{a}(db)-\rho(db)]\rho(da)\\
&=|\eta|^{-1}\sigma^{2}\sum_{1\leq i<j\leq m}\int_{[0,1]^{p}}[f(\textsl{x}^{i,j})-f(\textsl{x})]\rho^{\otimes m}(d\textsl{x}).
\end{align*}
We focus now on the immigration part. We take $f$ a function of the form $f:(x_{1},...,x_{m})\mapsto \prod_{i=1}^{m}\phi(x_{i})$ for some function $\phi$, and consider $F(\eta):= G_{f}(\rho)=\langle f, \rho^{\otimes m}\rangle$.
We may compute:
\begin{align}
F(\eta+h\delta_{a})-F(\eta)&=\left\langle \phi, \frac{\eta+h\delta_{a}}{z+h} \right\rangle^{m}-\langle \phi, \rho \rangle^{m} \nonumber \\
&=\sum_{j=2}^{m}\binom{m}{j}\left(\frac{z}{z+h}\right)^{m-j}\left(\frac{h}{z+h}\right)^{j}[\langle \phi,\rho\rangle^{m-j}\phi(a)^{j}-\langle \phi, \rho\rangle^{m}] \label{star} \\
&\qquad +m\left(\frac{z}{z+h}\right)^{m-1}\left(\frac{h}{z+h}\right)[\langle \phi,\rho\rangle^{m-1}\phi(a)-\langle \phi, \rho\rangle^{m}]\ \label{starstar}.
\end{align} We get that:
$$F'(\eta;a)=\frac{m}{z}\left[\phi(a)\langle\phi, \rho \rangle^{m-1}-\langle\phi, \rho \rangle^{m}\right].$$
Thus, 
\begin{equation*}
\label{factorize1}
F'(\eta;0)=|\eta|^{-1}\sum_{1\leq i\leq m}\int_{[0,1]^{p}}[f(\textsl{x}^{0,i})-f(\textsl{x})]\rho^{\otimes m}(d\textsl{x}) 
\end{equation*}
and
\begin{equation}
\int F'(\eta;a)\eta(da)=0 
\end{equation}
for such function $f$. 
This proves the Lemma.$\square$
\\
\\
This first lemma will allow us to prove the case $(i)$ of Theorem \ref{thm1}. We now focus on the case $(ii)$. Assuming that $\sigma^2=\beta=0$, the generator of $(M_{t}, t\geq 0)$ reduces to 
\begin{equation}
\label{discontinuous}
\mathcal{L}F(\eta)=\mathcal{L}_{0}F(\eta)+\mathcal{L}_{1}F(\eta)
\end{equation}
where, as in equations (3) and (4) of Subsection \ref{CBI},
\begin{align*} 
&\mathcal{L}_{0}F(\eta) = \int_{0}^{\infty}\hat{\nu_{0}}(dh)[F(\eta+h\delta_{0})-F(\eta)]\\ 
&\mathcal{L}_{1}F(\eta) = \int_{0}^{1}\eta(da)\int_{0}^{\infty}\hat{\nu_{1}}(dh)[F(\eta+h\delta_{a})-F(\eta)-hF'(\eta,a)].
\end{align*}
The following lemma is a first step to understand the infinitesimal evolution of the non-markovian process $(R_{t}, t\geq 0)$ in the purely discontinuous case.
\begin{lem} \label{genratio} Let $f$ be a continuous function on $[0,1]^{p}$ of the form $f(x_{1},...,x_{p})=\prod_{i=1}^{p}\phi(x_{i})$ and $F$ be the map $\eta \mapsto G_{f}(\rho)=\langle \phi,\rho \rangle^{p}$. Recall the notation $\rho:=\eta/ |\eta|$ and $z=|\eta|$. We have the identities: 
\begin{align*}
&\mathcal{L}_{0}F(\eta)=\int_{0}^{\infty}\hat{\nu_{0}}(dh)\left[G_{f}\left([1-\frac{h}{z+h}]\rho + \frac{h}{z+h}\delta_{0}\right)-G_{f}(\rho)\right] \\
&\mathcal{L}_{1}F(\eta)=z\int_{0}^{\infty}\hat{\nu_{1}}(dh)\int_{0}^{1}\rho(da)\left[G_{f}\left([1-\frac{h}{z+h}]\rho + \frac{h}{z+h}\delta_{a}\right)-G_{f}(\rho)\right].
\end{align*}
\end{lem}
\textit{Proof.}
The identity for $\mathcal{L}_{0}$ is plain, we thus focus on $\mathcal{L}_{1}$. 
Combining Equation (13) and the term (\ref{starstar}) we get
$$\int_{0}^{1} \rho(da)\left[ m\left(\frac{z}{z+h}\right)^{m-1}\left(\frac{h}{z+h}\right)[\langle \phi,\rho\rangle^{m-1}\phi(a)-\langle \phi, \rho\rangle^{m}]-hF'(\eta;a)\right]=0.$$ We easily check from the terms of \reff{star} that the map $h\mapsto \int_{0}^{1}\rho(da)[F(\eta+h\delta_{a})-F(\eta)-hF'(\eta,a)]$ is integrable with respect to the measure $\hat{\nu}_{1}$. This allows us to interchange the integrals and yields:
\begin{equation}
\label{factorize2}
\mathcal{L}_{1}F(\eta)=z\int_{0}^{\infty}\hat{\nu_{1}}(dh)\int_{0}^{1}\rho(da)\left[G_{f}\left(\frac{\eta+h\delta_{a}}{z+h}\right)-G_{f}(\rho)\right].
\end{equation}
$\square$\\
\\
The previous lemma leads us to study the images of the measures $\hat{\nu}_{0}$ and $\hat{\nu}_{1}$ by the map $\phi_{z}:h\mapsto r:=\frac{h}{h+z}$, for every $z>0$. Denote $\lambda_{z}^{0}(dr)=\hat{\nu}_{0}\circ \phi^{-1}_{z}$ and $\lambda_{z}^{1}(dr)=\hat{\nu}_{1}\circ \phi^{-1}_{z}$. 
The following lemma is lifted from Lemma 3.5 of \cite{Birk}. 
\begin{lem} \label{stablemes}
There exist two measures $\nu_{0}$, $\nu_{1}$ such that $\lambda_{z}^{0}(dr)=s_{0}(z)\nu_{0}(dr)$ and $\lambda_{z}^{1}(dr)=s_{1}(z)\nu_{1}(dr)$ for some maps $s_{0}, s_{1}$ from $\mathbb{R}_{+}$ to $\mathbb{R}$ if and only if for some $\alpha \in (0,2), \alpha' \in (0,1)$ and $c, c'>0$: $$\hat{\nu_{1}}(dx)=cx^{-1-\alpha}dx, \ \hat{\nu_{0}}(dx)=c'x^{-1-\alpha'}dx.$$ 
In this case: \begin{center} $s_{1}(z)=z^{-\alpha}$, $\nu_{1}(dr)=r^{-2}cBeta(2-\alpha,\alpha)(dr)$ \end{center} and \begin{center} $s_{0}(z)=z^{-\alpha'}$, $\nu_{0}(dr)=r^{-1}c'Beta(1-\alpha',\alpha')(dr)$.\end{center}
\end{lem}
\textit{Proof.} The necessary part is given by the same arguments as in Lemma 3.5 of \cite{Birk}. We focus on the sufficient part. Assuming that $\hat{\nu}_{0}, \hat{\nu}_{1}$ are as above, we have\\
\begin{itemize}
\item $\lambda^{1}_{z}(dr)=cz^{-\alpha}r^{-1-\alpha}(1-r)^{-1+\alpha}dr=z^{-\alpha} r^{-2}cBeta(2-\alpha, \alpha)(dr)$, and thus $s_{1}(z)=z^{-\alpha}.$
\item $\lambda_{z}^{0}(dr)=c'z^{-\alpha'}r^{-1-\alpha'}(1-r)^{-1+\alpha'}dr=z^{-\alpha'}r^{-1}c'Beta(1-\alpha',\alpha')(dr)$ and thus $s_{0}(z)=z^{-\alpha'}.$
$\square$
\end{itemize}
The next lemma allows us to deal with the second statement of Theorem \ref{thm1}. 
\begin{lem} \label{stable}
Assume that $\sigma^{2}=\beta=0$, $\hat{\nu}_{0}(dh)=ch^{-\alpha}1_{h>0}dh$ and $\hat{\nu}_{1}(dh)=ch^{-1-\alpha}1_{h>0}dh$. Let $f$ be a function on $[0,1]^{p}$ of the form $f(x_{1},...,x_{p})=\prod_{i=1}^{p}\phi(x_{i})$ , and $F$ be the map $\eta \mapsto G_{f}(\rho)$. We have 
$$|\eta|^{\alpha-1} \mathcal{L}F(\eta)=\mathcal{F}G_{f}(\rho),$$
for $\eta \neq 0$, where $\mathcal{F}$ is the generator of a $M$-Fleming-Viot process with immigration, with $M=(c'Beta(2-\alpha,\alpha-1),cBeta(2-\alpha,\alpha))$, see expressions $(3')$, $(4')$.
\end{lem}
\textit{Proof.} Recall Equation (\ref{discontinuous}):
\begin{equation*}
\mathcal{L}F(\eta)=\mathcal{L}_{0}F(\eta)+\mathcal{L}_{1}F(\eta)
\end{equation*}
Recall from Equation (13) that we have $\int_{0}^{1}F'(\eta;a)\eta(da)=0$ for $F(\eta)=G_{f}(\rho)$. Applying Lemma \ref{genratio} and Lemma \ref{stablemes}, we get that in the case $\sigma^2=\beta=0$ and $\hat{\nu_{1}}(dx)=cx^{-1-\alpha}dx, \hat{\nu_{0}}(dx)=c'x^{-1-\alpha'}dx$:  
\begin{align*}
\mathcal{L}F(\eta)=\mathcal{L}G_{f}(\rho)=& \quad s_{0}(z)\int_{0}^{1}r^{-1}c'Beta(1-\alpha',\alpha')(dr)[G_{f}((1-r)\rho+r\delta_{0})-G_{f}(\rho)]\\
&+zs_{1}(z)\int_{0}^{1}r^{-2}cBeta(2-\alpha,\alpha)(dr)\int_{0}^{1}\rho(da)[G_{f}((1-r)\rho+r\delta_{a})-G_{f}(\rho)].
\end{align*} 
Recalling the expressions (3'), (4'), the factorization $h(z)\mathcal{L}F(\eta)=\mathcal{F}G(\rho)$ holds for some function $h$ if $$s_{0}(z)=zs_{1}(z),$$ if $\alpha'=\alpha-1$. In that case, $h(z)=z^{\alpha-1}$. $\square$
\\
\\
We are now ready to prove Theorem \ref{thm1}. To treat the case (i), replace $\alpha$ by $2$ in the sequel. The process $(Y_{t}, R_{t})_{t\geq 0}$ with lifetime $\tau$ has the Markov property. The additive functional $C(t)=\int_{0}^{t}\frac{1}{Y_{s}^{\alpha-1}}ds$ maps $[0,\tau)$ to $[0,\infty)$. From Theorem 65.9 of \cite{Sharpe} and Proposition \ref{timechange}, the process $(Y_{C^{-1}(t)}, R_{C^{-1}(t)})_{t\geq 0}$ is a strong Markov process with infinite lifetime. Denote by $\mathcal{U}$ the generator of $(Y_{t}, R_{t})_{t\geq 0}$. As explained in Birkner et al. \cite{Birk} (Equation (2.6) p314), the law of $(Y_{t},R_{t})_{t\geq 0}$ is characterized by $\mathcal{U}$ acting on the following class of test functions: \begin{center} $(z, \rho)\in \mathbb{R}_{+}\times\mathcal{M}_{1} \mapsto F(z, \rho):=\psi(z)\langle\phi, \rho\rangle^{m}$ \end{center} for $\phi$ a non-negative measurable bounded function on $[0,1]$, $m\geq 1$ and $\psi$ a twice differentiable non-negative map. Theorem 3 of Volkonski{\u\i}, see \cite{timechange} (or Theorem 1.4 Chapter 6 of \cite{EthierKurtz}) states that the Markov process with generator $$\tilde{\mathcal{U}}F(z,\rho):=z^{\alpha-1}\mathcal{U}F(z,\rho)$$ 
coincides with $(Y_{C^{-1}(t)},R_{C^{-1}(t)})_{t\geq 0}$.
We establish now that $(R_{C^{-1}(t)}, t\geq 0)$ is a Markov process with the same generator as the Fleming-Viot processes involved in Theorem \ref{thm1}. Let $G(z,\rho)=G_{f}(\rho)=\langle\phi, \rho\rangle^{m}$ (taking $f:(x_{1},...,x_{m})\mapsto \prod_{i=1}^{m}\phi(x_{i})$). In both cases (i) and (ii) of Theorem \ref{thm1}, we have:
\begin{align*}
z^{\alpha-1}\mathcal{U}G(z,\rho)&=z^{\alpha-1}\mathcal{L}F(\eta) \text{ with } F: \eta \mapsto G_{f}(\rho)\\
&=\mathcal{F}G_{f}(\rho).
\end{align*}
First equality holds since we took $\psi\equiv 1$ and the second uses Lemma \ref{continuous} and Lemma \ref{stable}.
Since it does not depend on $z$, the process $(R_{C^{-1}(t)}, t\geq 0)$ is a Markov process, moreover it is a generalized Fleming-Viot process with immigration with parameters as stated. $\square$
\\
\\
\textit{Proof of Proposition \ref{independence}.} Let $(Y_{t})_{t\geq 0}$ be a Feller branching diffusion with continuous immigration with parameters $(\sigma^{2}, \beta)$. Consider an independent $M$-Fleming-Viot $(\rho_{t}, t\geq 0)$ with $M=(\beta \delta_{0}, \sigma^{2}\delta_{0})$. We first establish that $(Y_{t}\rho_{C(t)}, 0\leq t< \tau)$ has the same law as the measure-valued branching process $(M_{t}, 0\leq t< \tau)$.
Recall that $\mathcal{L}$ denote the generator of $(M_{t}, t\geq 0)$ (here only the terms (1) and (2) are considered). Consider $F(\eta):=\psi(z)\langle \phi, \rho\rangle ^{m}$ with $z=|\eta|$, $\psi$ a twice differentiable map valued in $\mathbb{R}_{+}$ and $\phi$ a non-negative bounded measurable function. Note that the generator acting on such functions $F$ characterizes the law of $(M_{t\wedge \tau}, t\geq 0)$.
First we easily obtain that
\begin{align*}
F'(\eta;0)&=\psi'(z)\langle \phi, \rho \rangle^{m}+m  \frac{\psi(z)}{z}[\phi(0)\langle \phi, \rho \rangle^{m-1}-\langle \phi, \rho \rangle^{m}],\\
F''(\eta; a,b)&=\psi''(z)\langle \phi, \rho \rangle^{m}+m\frac{\psi'(z)}{z}\left[(\phi(b)+\phi(a)) \langle \phi, \rho \rangle^{m-1}-2\langle \phi, \rho \rangle^{m} \right] \\
&\qquad+m \frac{ \psi(z) }{z^{2}}   \left[ (m-1)\phi(a)\phi(b)\langle \phi, \rho\rangle^{m-2} - m (\phi(a)+\phi(b)) \langle \phi, \rho\rangle ^{m-1} + (m+1) \langle \phi, \rho\rangle ^{m}\right] .
\end{align*} 
Simple calculations yield, \begin{align*} \mathcal{L}F(\eta)&=\left[z\left(\frac{\sigma^{2}}{2}\psi''(z)\right)+\beta\psi'(z)\right]\langle \phi, \rho\rangle^{m}\\
&\qquad+\frac{\psi(z)}{z}\left[\sigma^{2}\frac{m(m-1)}{2}\left(\langle\phi^{2},\rho\rangle\langle \phi, \rho\rangle^{m-2}-\langle \phi, \rho \rangle^{m}\right)+ \beta m \left(\phi(0)\langle \phi,\rho \rangle^{m-1}-\langle \phi,\rho \rangle^{m} \right) \right].
\end{align*}
We recognize in the first line the generator of $(Y_{t}, t\geq 0)$ and in the second, $\frac{1}{z}\mathcal{F}G_{f}(\rho)$ with $f(x_{1},...,x_{m})=\prod_{i=1}^{m}\phi(x_{i})$ and $c_{0}=\beta, c_{1}=\sigma^{2}$. We easily get that this is the generator of the Markov process $(Y_{t}\rho_{C(t)}, t\geq 0)$ with lifetime $\tau$. 
We conclude that it has the same law as $(M_{t\wedge \tau}, t\geq 0)$. We rewrite this equality in law as follows:
\begin{equation}
\label{decomposition}
(Y_{t} \;  \rho_{C(t)}, 0\leq t<\tau)  \overset{law}{=}  (|M_{t}|\; R_{C^{-1}(C(t))}, 0 \leq t < \tau), 
\end{equation}
with $C$ defined by $C(t)=\int_{0}^{t} |M_{s}|^{-1} ds$ for $0 \leq t < \tau$ on the right hand side. Since $(C(t), t\geq 0)$ and $(\rho_{t}, t\geq 0)$ are independent on the left hand side and the decomposition in (\ref{decomposition}) is unique, we have also $(C(t), 0\leq t< \tau)$ and $(R_{C^{-1}(t)}, 0\leq t<\tau)$ independent on the right hand side. \\
\\
Concerning the case $(ii)$ of Theorem \ref{thm1}, we easily observe that the presence of jumps implies that such a decomposition of the generator cannot hold. See for instance Equation (2.7) of \cite{Birk} p344. The processes $(R_{C^{-1}(t)}, t\geq 0)$ and $(Y_{t}, \geq 0)$ are not independent.
$\square$
\\
\\
\textbf{Acknowledgments}. The authors would like to thank Jean Bertoin and Jean-Fran\c cois Delmas for their helpful comments and advice. C.F thanks the Statistical Laboratory of Cambridge where part of this work was done with the support of the Foundation Sciences Math\'ematiques de Paris.  O.H. thanks Goethe Universit\"at of Frankfurt for hospitality, and \'Ecole Doctorale MSTIC for support. This work is partially supported by the ``Agence Nationale de la Recherche'', ANR-08-BLAN-0190.

\end{document}